\documentclass[12pt,fleqn]{article}
\usepackage{amsfonts}
\usepackage{amsmath}
\usepackage{amssymb}

\newcommand{\con}{\mathfrak c}

\newcommand{\eps}{\varepsilon}

\newcommand{\algb}{\mathfrak B}

\newcommand{\BB}{\protect{\mathcal B}}
\newcommand{\AAA}{\protect{\cal A}}
\newcommand{\CC}{{\cal C}}
\newcommand{\FF}{\protect{\mathcal F}}

\newcommand{\UU}{{\cal U}}
\newcommand{\JJ}{{\cal J}}
\newcommand{\HH}{{\cal H}}

\newcommand{\LL}{\protect{\mathcal D}}

\newcommand{\btu}{\bigtriangleup}

\newcommand{\dist}{\protect{\rm dist}}

\newcommand{\fin}{\protect{\rm Fin}}

\newcommand{\er}{\mathbb R}
\newcommand{\qu}{\mathbb Q}

\newcommand{\vf}{\varphi}
\newcommand{\osc}{{\rm osc}}
\newcommand{\ol}{\overline}

\newcommand{\sm}{\setminus}
\newcommand{\sub}{\subseteq}

\newcommand{\stevo}{Todor\v{c}evi\'c}
\newenvironment{mypr}{\noindent {\em Proof.}\ }{\noindent $\diamondsuit$\medskip}

\newcommand{\conv}{\operatorname{conv}}

\newtheorem{thm}{Theorem}[section]

\newtheorem{lem}[thm]{Lemma}
\newtheorem{cor}[thm]{Corollary}

\newenvironment{ex}
{
\bigskip
\addtocounter{thm}{1}
\noindent
{\bf Example \arabic{section}.\arabic{thm}}
}
{$\diamondsuit$
\bigskip}

\newenvironment{remark}
{
\bigskip
\addtocounter{thm}{1}
\noindent
{\bf Remark \arabic{section}.\arabic{thm}}
}
{$\diamondsuit$
\bigskip}

\newcounter{nrpun}
\newenvironment{liro}
{
\begin{list}{\it (\roman{nrpun})}{\usecounter{nrpun}
\leftmargin1.4cm \rightmargin0cm \topsep0.3cm\itemsep0pt}
}
{\end{list}}

\newenvironment{liar}
{
\begin{list}{\it (\arabic{nrpun})}{\usecounter{nrpun}
\leftmargin1.4cm \rightmargin0cm \topsep0.3cm\itemsep0pt}
}
{\end{list}}
\newenvironment{mysec}{\medskip\pagebreak[2]
\addtocounter{section}{1}\sc \begin{center}\arabic{section}. }
{\setcounter{thm}{0}\nopagebreak \end{center} \medskip}

\date{}

\title{On measures on Rosenthal compacta }
\author{Witold Marciszewski and Grzegorz Plebanek}
\date{}

\begin{document}

\maketitle

\renewcommand{\thefootnote}{\fnsymbol{footnote}\noindent}
\footnotetext{2000 {\em Mathematics Subject Classification}:
 Primary 28C15, 46A50; Secondary 28A60, 54C35.}
\footnotetext{{\em Key words and phrases:} Radon measure, countably determined measure, Rosenthal compact.}
\footnotetext{Research of the first author was partially supported by
MNiSW Grant Nr N N201 382034.}
\footnotetext{The second author was partially supported by the grant 2191/W/IM/09
from University of Wroc\l aw.}

\begin{abstract}
We show that if $K$ is Rosenthal compact which can be represented
by functions with countably many discontinuities then every Radon
measure on $K$ is countably determined. We also present an
alternative proof of the result stating that every Radon measure
on an arbitrary Rosenthal compactum is of countable type. Our
approach is based on some caliber-type properties of measures,
parameterized by separable metrizable spaces.
\end{abstract}

\begin{mysec}
Introduction
\end{mysec}

For any compact space $K$ we denote by $P(K)$ the space
of Radon probability measures on $K$.
Recall that a measure $\mu\in P(K)$ is {\em countably determined} if
there is a countable family $\FF$ of closed sets such that
\[\mu(U)=\sup\{ \mu(F): F\sub U, F\in\FF\},\]
for every open $U\sub K$; such a family $\FF$ is said
to approximate from below all open sets. Countably determined measures were considered by
Pol \cite{Po82} and Mercourakis \cite{Me96}.

It is clear that a countably determined measure is of countable Maharam type. Recall that a measure
$\mu$ is said to have countable (Maharam) type if $L_1(\mu)$ is separable, which
is equivalent to saying  that there is a countable family $\FF\sub Bor(K)$ which is $\btu$-dense
in $Bor(K)$.
While the Maharam type of a measure is a cardinal coefficient of its measure algebra, the property
of being countably determined describes a connection between the measure and the topology of the underlying
space.
For instance the Lebesgue measure $\lambda$ on $[0,1]$ is  countably determined but the
 measure $\widehat{\lambda}$ that can be defined
on the Stone space of the measure algebra of $\lambda$ is not countably determined though still
of countable type.

Throughout this note $X$ denotes a separable metrizable space, and
$B_1(X)$ stands for the space of first  Baire class functions $X\to
\er$, equipped with the pointwise topology. A compact space $K$ is
said to be {\em Rosenthal compact} if $K$ can be topologically
embedded into $B_1(X)$ for some Polish space $X$; cf.\ Godefroy
\cite{Go80}, Bourgain et al.\ \cite{BFT}, Marciszewski
\cite{Ma92}, \stevo\ \cite{To99}.

We address here the following question asked by Pol \cite{Po82}: suppose that $K$ is
Rosenthal compact; is every $\mu\in P(K)$ countably determined?
Mercourakis \cite{Me96} mentions several
classes of compact spaces on which every regular measure is countably determined and
also asks if all Rosenthal compacta have such a property. In  section 4 we present
a partial affirmative answer to Pol's question, see in particular Theorem \ref{4:2}.

Let us recall other measure-theoretic properties of Rosenthal compacta that have been
considered.
If $K$ is Rosenthal compact then  by a result due to Godefroy \cite{Go80} the space $P(K)$ is
again Rosenthal compact in its $weak^*$ topology.  This fact, together with
the Frechet property of Rosenthal compacta (see \cite{BFT}), imply that
the support of each $\mu\in P(K)$ is separable. Talagand mentions without proof
in \cite{Ta84}, (14-2-2) the following two unpublished results of Bourgain \cite{Bo76}:

\begin{liar}
\item \label{1a} every $\mu\in P(K)$ is of countable Maharam type;
\item \label{1b} if $K\sub B_1(X)$ for some Polish space $X$ then
for any $\mu\in P(K)$ the mapping $e: X\to L_1(\mu)$,
$e(x)(f)=f(x)$, for $f\in K, x\in X$, is of the first Baire class.
\end{liar}

In fact, in his thesis Bourgain gave a detailed proof of (\ref{1b}), and later only
announced (\ref{1a}),  with a comment  that it
can be obtained by a modification of his argument leading to (\ref{1b}), see
\cite[p. 33]{Bo76}.\footnote{We wish to thank Antonio
Avil\'{e}s for providing us with a copy of \cite{Bo76}}

Stevo \stevo\ presented an interesting argument for (\ref{1a}) using a result due to Fremlin
 \cite{Fr97} that under MA$(\omega_1)$ every compact space carrying a Radon measure of uncountable
 Maharam type can be continuously mapped onto $[0,1]^{\omega_1}$. Since
no Rosenthal compactum admits such a surjection, it follows from Fremlin's result that
the sentence
 `every Radon measure on a Rosenthal compactum has countable type' is relatively consistent.
 \stevo\ in \cite{To99} (see also \cite{To00}) analyzed properties of Rosenthal compacta
 preserved in forcing extensions and building on this was able to give a proof of (\ref{1a}) requiring
no extra set-theoretic assumptions.

In Theorem \ref{4:3} we present a somewhat more direct proof of
Bourgain's result (\ref{1a}). Our approach is based on some
combinatorial results related to Polish spaces given in section 3
below; we also use some result related to measure algebras, stated
as Theorem \ref{3:3} , which originated in Fremlin's paper
\cite{Fr97} and was later generalized by Fremlin and Plebanek in
an unpublished preprint \cite{FP04}. The self-contained proof of
Theorem \ref{3:3} is enclosed in section 5 for completeness.

We would like to comment on similarities and differences between
our proof of (1) and \stevo's argument from \cite{To99}. While we
use at times widely accepted terminology connected with forcing,
our approach is entirely based on some auxiliary combinatorial
results on measure and topology, and requires no considerations
related to properties preserved in forcing extensions. A
mathematician with an expertise in forcing may feel that such a
difference is not very essential but we hope that a substantial
mathematical audience will find our proof useful.

\begin{mysec}
Measures on function spaces
\end{mysec}
Given a function $g:X\to\er$, $Y\sub X$ and $x\in Y$, we write
$\osc(g,Y,x)$ for the oscillation of the restricted function
$g_{|Y}$ at the point $x$ (and $\osc(g,x)$ in case $Y=X$).
Observe that for a separable metrizable space $X$ and a subspace
$K$ of the product $\er^X$, for every $Y\sub X, x\in Y$, and
$\delta>0$, the set
\[\{g\in K:\; \osc(g,Y,x)\ge \delta\},\]
 is a
$G_\delta$-subset of $K$.

We shall need the following remark on countable determinacy of measures:
if $\HH$ is a pseudobase of the topology on $K$, and $\FF$ is
a countable family approximating every $U\in\HH$ from below with respect to
a fixed measure $\mu\in P(K)$ then $\mu$ is countably determined. Indeed,
it is routine to check that the family $\LL=(\FF)_{\cap\,\cup}$ (i.e.\ the lattice generated by $\FF$)
approximates from below
every element from the lattice $\HH'$ generated by $\HH$.
In turn $\HH'$ is a base for the topology that is closed under finite unions, so
$\HH'$ approximates from below all open sets by regularity of $\mu$.
As $\LL$ is countable, this shows that $\mu$ is countably determined.

\begin{lem}\label{2:1}
Let $X$ be a separable metrizable space, $K\sub \er^X$ be a
compact space and $\mu\in P(K)$. Suppose that for every $\delta>0$
the space $X$ can be written as
$X=\bigcup_{n\in\omega}X_n(\delta)$ so that for every $n$ and
$x\in X_n(\delta)$ we have
\[\mu\{g\in K:\; \osc(g,X_n(\delta),x)\ge \delta\}<\delta.\]
Then the measure $\mu$ is countably determined.
\end{lem}

\begin{mypr} Let us fix a countable base $\UU$ of $X$ and denote by $\JJ$ the family of all closed intervals
with rational endpoints. For any $I\in\JJ$ and  $E\sub X$ we write
\[A(I,E)=\{g\in K:\; g[E]\sub I\};\]
note that such a set $A(I,E)$ is closed in $K$ for arbitrary $E\sub X$.

We apply the assumption of the lemma to every $\delta$ from $N=\{1,1/2,1/3,\ldots\}$, and consider
the family $\FF$ of all sets of the form
\[A(I,X_n(\delta)\cap U)\mbox{ where } I\in\JJ, n\in\omega, \delta\in N, U\in\UU.\]
We shall check that the lattice $\LL$ generated by $\FF$, which is clearly countable, approximates
all open sets in $K$ from below.
By the remark preceding Lemma \ref{2:1} it will do to approximate from below
an open subset $H\sub K$ of the form
$H=\{g\in K:\; g(x)\in J\}$, for a fixed $x\in X$ and an open interval $J=(t,s)\sub\er$.

Let $\eps>0$. First find $\delta>0$ such that writing
\[J'=(t+2\delta,s-2\delta) \mbox{ and } H'=\{g\in K: g(x)\in J'\},\]
 we have $\mu(H')>\mu(H)-\eps/3$. We may also assume that $\delta<\eps/3$ and $\delta\in N$.

Choose $n$ such that $x\in X_n(\delta)$ and write $M=\{g\in K: \osc(f,X_n(\delta),x)< \delta\}$.
Then $\mu(M)\ge 1-\delta$ by the assumption of the lemma so
\[\mu(H\cap M)>\mu(H)-\delta \mbox{ and } \mu(H'\cap M)>\mu(H\cap M)-\eps/3> \mu(H)- 2\eps/3.\]
Take any rational numbers $p\in (t, t+\delta)$ and $q\in (s-\delta, s)$, and let $I=[p,q]$.
Observe that whenever $g\in H'\cap M$ then there is
$U\in \UU$ containing $x$ such that $g[U\cap X_n(\delta)]\sub I$, so $g\in A(I,U\cap X_n(\delta))\sub H$.

It follows that there are $A_n\in \LL$ such that $H'\cap M\sub \bigcup_{n\in\omega}A_n\sub H$. Then
$B=\bigcup_{n\le k} A_n\in \LL$ and
for $k$ sufficiently large we have
\[\mu(B)> \mu(H'\cap M)-\eps/3>  \mu(H)-\eps,\]
 and this completes the proof.\end{mypr}

\begin{cor}\label{2:2}
Let $X$ be a separable metrizable space, $K\sub \er^X$ be a
compact space and $\mu\in P(K)$. Suppose that there is a
decomposition $X=\bigcup_n X_n$ such that for every $n$ and every
$x\in X_n$ the function $g_{|X_n}$ is continuous at $x$ for
$\mu$-almost all $g\in K$. Then the measure $\mu$ is countably
determined.
\end{cor}

We do not know if the assumption of Lemma \ref{2:1} is fulfilled
by every measure on a Rosenthal compact space $K$. The following
observation might be useful when analyzing that.

\begin{lem}\label{2:3}
Let $X$ be a metric space, $K\sub \er^X$ be a compact space and $\mu\in P(K)$.
Then for any $\eps,\delta>0$ the set
\[X_0=\{x\in X: \mu(\{g\in K: \osc(g,x)\ge\delta\})\ge\eps,\]
is closed.
\end{lem}

\begin{mypr}
Take $x_n\in X_0$ and suppose that $x_n\to x\in X$. Write
\[B_n=\{g\in K: \osc(g,x_n)\ge\delta\},\quad B=\bigcap_{n}\bigcup_{k\ge n} A_k.\]
Then $\mu(B_n)\ge\eps$ for every $n$ so $\mu(B)\ge\eps$ as well;
it suffices to notice that for every $g\in B$ we have
$\osc(g,x)\ge\delta$.
\end{mypr}

We shall show that Corollary \ref{2:2} is applicable
to those  $K$ that can be represented by
functions with few points of discontinuity, see Theorem \ref{4:2}. For that purpose we need some
combinatorial results given in the next section.

\begin{mysec}
Parameterized calibers
\end{mysec}

In this section we consider a fixed compact space $K$ and a
probability Borel measure $\mu$ on $K$. A family $\{L_t: t\in T\}$
of (measurable) subsets of $K$ is said to be {\em $\mu$-centered}
if
\[\mu(\bigcap_{t\in a} L_t)>0,\]
 for every finite set $a\sub T$. Here $T$ was chosen to denote some index set
but in the sequel  $T$ will often denote  a separable metrizable
space. In such a case, we denote by $\fin(T)$ the space of all finite
subsets of $T$ endowed with the Vietoris topology. Recall that a
basic open  set in $\fin(T)$ is of the form
\[\{a\in \fin(T): a\sub \bigcup_{i\le n} V_i,\; a\cap V_i\neq\emptyset \mbox{ for all } i\le n\},\]
where the sets $V_i$ are open in $T$.

Suppose now that $T\ni t\to L_t\sub K$ is an arbitrary mapping,
where $\mu(L_t)>0$ for every $t\in T$. It is not difficult to
check that if $T$ is uncountable then there is an infinite set
$S\sub T$ such that $\{L_t: t\in S\}$ is $\mu$-centered. It is
also well-known that the existence of an uncountable $S$ with such
a property is undecidable within the usual axioms of set theory,
see e.\ g.\ D\v{z}amonja \& Plebanek \cite{DP04}. We shall
investigate here  if  sets $S\sub T$ giving rise to centered
families can have some topological properties as subspaces of $T$.

A family $\AAA\sub \fin(T)$ will be called a $ccc$-family if for every uncountable pairwise
disjoint subfamily $\AAA_0\sub\AAA$ there are distinct $a,b\in\AAA_0$ such that $a\cup b\in\AAA$.

\begin{lem}\label{3:1}
Let $T$ be an uncountable separable metrizable space and let
$\AAA\sub \fin(T)$ be a hereditary $ccc$-family such that
$\bigcup\AAA=T$. Then there is a countable dense-in-itself set
$S\sub T$ such that  $\fin(S)\sub\AAA$.
\end{lem}

\begin{mypr}
The proof is based on the following observation.

\noindent
{\bf (1) Claim.} There is a countable set $J\sub T$ such that every $a\in\AAA$ with
$a\cap J=\emptyset$ has the following property (*):
\medskip

(*)\quad For every open neighborhood $U$ of $a\in\fin(T)$ there is
an uncountable pairwise disjoint family $\BB\sub \AAA\cap U$ such
that $a\cup b\in \AAA$ for every $b\in \BB$.
\medskip

Suppose converse; then by an obvious transfinite induction we can find a pairwise disjoint
uncountable family $\BB\sub\AAA$ such
that no $a\in\BB$ satisfies (*). Let $\UU$ be a countable base of $\fin(T)$. For every $a\in\BB$ choose
$U_a\in\UU$ witnessing that $a$ does not satisfy (*). As $\UU$ is countable, shrinking
$\BB$ if necessary,  we can assume that for all $a\in\BB$, $U_a=U$ with $U$ fixed. Note that for a given
$a\in\BB$ the family
\[\{b\in\BB: b\sub U, b\cup a\in \AAA\},\]
is countable. Therefore we can choose an uncountable family $\CC\sub\BB$ such that for each pair
$a,b$ of distinct elements of $\CC$ we have $a\cup b\notin \AAA$, a contradiction.
\medskip

\noindent
{\bf (2)} Using (1) we can construct inductively  sets $a_n\in\AAA$ such that for all $n$

\begin{liro}
\item $a_n\cap J=\emptyset$;
\item $a_n\sub a_{n+1}$
\item $a_n$ is $1/n$--dense-in-itself, i.e.\ for every $p\in a_n$ there is $q\in a_n$ such that
$q\neq p$ and $\dist(p,q)<1/n$.
\end{liro}

Now the set $S=\bigcup_{n\in\omega} a_n$ is dense-in-itself and
all finite subsets of $S$ are in $\AAA$, as required.
\end{mypr}

\begin{thm}\label{3:2}
If $T\ni t\to L_t\sub K$ is any mapping from an uncountable
separable metrizable space $T$ into the family of compact sets of
positive measure then $\bigcap_{t\in S} L_t\neq\emptyset$ for some
dense-in-itself set $S\sub T$.
\end{thm}

\begin{mypr}
For any $a\in\fin(T)$ we  denote $L_a=\bigcap_{t\in a} L_t$ and set
\[\AAA=\{a\in \fin(T): \mu(L_a)>0\}.\]
Then $\AAA$ is a hereditary family and $\{t\}\in\AAA$ for every
$t\in T$. If $\AAA_0\sub\AAA$ is any uncountable subfamily then
$\{L_a: a\in\AAA_0\}$ is an uncountable family of nonnull sets so
$\mu(L_a\cap L_b)>0$ for some distinct $a,b\in\AAA_0$; then $a\cup
b\in\AAA$. Hence $\AAA$ is a $ccc$-family and by Lemma \ref{3:1}
there is a dense-in-itself set $S\sub T$ such that $\{L_s: s\in
S\}$ is $\mu$-centered, and $\bigcap_{t\in S} L_t\neq\emptyset$ by
compactness.
\end{mypr}

Let us observe that it is easy to prove Theorem \ref{3:2} assuming
MA + non CH: indeed, then $\omega_1$ is a precaliber of measures,
i.e.\  there is uncountable $T_0\sub T$ such that the family
$\{L_t: t\in T_0\}$ is $\mu$-centered;  the assertion of \ref{3:2}
follows from the fact that $T_0$ may have only countably many
isolated points. A ZFC proof presented above is a modification of
an argument from \cite{Ma92}, Theorem 4.1.

The following example shows that the assertion of Theorem
\ref{3:2} cannot be strengthened by replacing "dense-in-itself"
with "not nowhere dense". Let us point out that if such a stronger
result were true then one could show that every regular
probability measure on an arbitrary Rosenthal compactum is
countably determined, cf.\ the proof of Theorem \ref{4:2}.

\begin{ex} In the setting of Theorem \ref{3:2}, it may happen that every set $S\sub T$ such that
$\bigcap_{t\in S} L_t\neq\emptyset$  is necessarily nowhere dense.

We let $T=2^\omega$ and consider the standard product measure $\mu$ on $K=2^\omega$. Let us fix a sequence
$(I_n)_{n\ge 1}$ of pairwise disjoint subset of $\omega$, such that $|I_n|=n+1$ for $n=1,2,\ldots$.
For any finite $I\sub\omega$ and $\vf:I\to 2$ we write
\[C(I,\vf)=\{t\in 2^\omega: t(i)=\vf(i)\mbox{ for } i\in I\},\]
for the corresponding cylinder set; note that
\[ \mu(C(I,\vf))=1/2^{|I|}.\]
For $t\in 2^\omega$ we define $L_t$ by the formula
\[L_t=2^\omega\sm \bigcup_{n\ge 1} C(I_n,t_{|I_n}\};\]
simple calculations show that $\mu(L_t)\ge 1/2$ for every $t$.

Suppose that $S\sub 2^\omega$ is not nowhere dense, i.e.\ there is
some basic open set $C(I,\vf)$ which is contained in the closure
of $S$. Take any $t\in 2^\omega$ and $n$ such that $I\cap
I_n=\emptyset$. Then the set
\[S'=S\cap C(I,\vf)\cap C((I_n,t_{|I_n})\neq\emptyset;\]
if $s\in S'$ then $s_{|I}=\vf$ and $s_{|I_n}=t_{|I_n}$, which
implies $t\notin L_s$. In this way  we have checked that
$\bigcap_{s\in S}L_s=\emptyset$.
\end{ex}

Let us consider now an indexed family $\{(L_t^0,L_t^1): t\in T\}$ of disjoint pairs of sets
$L_t^0,L_t^1\sub K$. Such a family is called {\em independent}
if for every finite set $a\sub T$ and every function $\vf:a\to 2$
\[L_\vf=\bigcap_{t\in a} L_t^{\vf(t)}\neq\emptyset.\]
In a similar way we define  the $\mu$-independence, which means that for some measure $\mu$ on $K$
every set $L_\vf$ as above is rather of positive measure than simply nonempty.

We shall need the following technical result.

\begin{thm}[\cite{FP04}]\label{3:3}
Let $\mu$ be a probability measure on a space $K$.
Suppose that
$\{(L^0_\xi,L^1_\xi):\xi<\omega_1\}$ is a family of disjoint pairs of  measurable subsets of $K$ such
that for some constant $\eps>0$

\begin{liro}
\item $\mu(L_\xi^0)+\mu(L_\xi^1)>1-\eps/2$ for every $\xi<\omega_1$;
\item $ \mu(L^0_\xi\cap L^1_\eta)\ge \eps$ whenever $\xi,\eta<\omega_1$, $\xi\neq \eta$.
\end{liro}

Let $\AAA$ be a family of those finite sets $a\sub \omega_1$ for
which $\{(L_\xi^0,L_\xi^1): \xi\in a\}$ is $\mu$-independent. Then
there are an uncountable set $T\sub\omega_1$, and a hereditary
$ccc$-family $\AAA_0\sub \AAA\cap {\mathcal P}(T)$ such that
$\bigcup\AAA_0=T$.
\end{thm}

A special case of Theorem \ref{3:3}  appeared in
Fremlin \cite{Fr97} (see the proof of Theorem 6, there); in
its present form the result can be derived
from the argument given in an unpublished note by Fremlin \& Plebanek \cite{FP04}.
We enclose a self-contained proof of \ref{3:3} in the last section.

\begin{cor} \label{3:4}
Let $T$ be an uncountable separable metrizable space and let
$\mu$ be a probability measure on a space $K$. Suppose that,
$\{(L^0_t,L^1_t):t\in T\}$ is a family of disjoint pairs of
measurable subsets of $K$ such that for some constant $\eps>0$

\begin{liro}
\item $\mu(L_t^0)+\mu(L_t^1)>1-\eps/2$ for every $t\in T$;
\item $ \mu(L^0_t\cap L^1_s)> \eps$ whenever $t,s\in T$, $t\neq s$.
\end{liro}

Then there is a dense-in-itself set $S\sub T$ such that
the family $\{(L_t^0,L_t^1):t\in S\}$ is independent.
\end{cor}

\begin{mypr}
We apply Theorem \ref{3:3} to get an uncountable $T_0\sub T$ and a $ccc$-family
$\AAA_0\sub  \fin(T_0)$ such that $\bigcup\AAA_0=T_0$ and
$\{(L_t^0,L_t^1):t\in a\}$ is independent for every $a\in\AAA_0$.
Now the assertion of the theorem follows immediately from Lemma \ref{3:1}.
\end{mypr}

\begin{mysec}
Applications to Rosenthal compacta
\end{mysec}

For a (separable metrizable) space $X$ we write $CD(X)$ for the
space of all functions $g:X\to \er$ for which the set of points of
discontinuity is at most countable. The following fact is due to
Marciszewski and Pol \cite[Proposition 2.2]{MP09}.

\begin{thm}\label{4:1}
Let $X$ be a Borel subspace of a separable completely
 metrizable space.
Every compact space $K\sub CD(X)$ can be embedded into
$CD(2^\omega)$.
\end{thm}

\begin{thm} \label{4:2}
Let $X$ be a Borel subspace of a separable completely
 metrizable space and suppose that $K\sub CD(X)$ is a compact space.
Then for every $\mu\in P(K)$ the set
\[\{x\in X:\mu\{g\in K:\; \osc(g,x)>0\}>0\},\]
is countable. Consequently, every measure $\mu\in P(K)$ is countably determined.
\end{thm}

\begin{mypr}
By Theorem \ref{4:1} we can assume that $K\sub CD(2^\omega)$. Suppose that $X'\sub 2^\omega$ is
uncountable and
\[\mu\{g\in K:\; \osc(g,x)>0\}>0,\]
 for every $x\in X'$. Then there is uncountable
$T\sub X'$ and $\eps>0$ such that the set
\[D_t=\{g\in K:\; \osc(g,t)\ge \eps\},\]
has positive measure for $t\in T$. By regularity of $\mu$, for
every $t\in T$ there is a compact set $L_t\sub D_t$ with
$\mu(L_t)>0$.

By Theorem \ref{3:2} there is a dense-in-itself set $S\sub 2^\omega$ such that
\[L_S=\bigcap_{x\in S} L_x\neq\emptyset.\]
 But if $g\in L_S$ then $g$ is clearly discontinuous at each
$x\in \overline{S}$; since $\overline{S}$ is a perfect subset of
the Cantor set it has size $\con$, a contradiction. The second assertion follows from Corollary \ref{2:2}.
\end{mypr}

\begin{remark} Pol \cite{Po82} and Mercourakis \cite{Me96} considered another
property of measures:  $\mu\in P(K)$ is {\em strongly countably determined} if there is
a countable family of closed $G_\delta$ subsets of $K$, approximating all open sets in $K$ from below
(with respect to $\mu$). Note that the Dirac measure $\delta_x$, where $x\in K$, is always countably determined
while $\delta_x$ is strongly countably determined if and only if $x$ is a $G_\delta$ point in $K$.
Strongly countably determined measures were introduced by Babiker \cite{Ba76} (under the name {\em uniformly
regular measures}); see Plebanek \cite{Pl00} for further results and references.

Let $H$ be the Helly space of all
nondecreasing functions from $[0,1]$ into $[0,1]$. Observe, that
for this compact space, the sets $A(I,E)$ defined in the proof of
Lemma \ref{2:1} are $G_\delta$-sets in $H$. Indeed, for any
$E\subset [0,1]$ we can find a countable subset $F$ of $E$ such
that $\conv F = \conv E$. Then $A(I,E) = A(I,F)$, since all
functions form $H$ are nondecreasing.  Therefore, proofs of Lemma
\ref{2:1} and Theorem \ref{4:2} show that every Radon probability
measure on $H$ is strongly countably determined.

Let $BV$ denote another well-known example of a Rosenthal
compactum, the space of all functions from $[0,1]$ into $[0,1]$ of
total variation $\le 1$. Every function $f\in BV$ can be
represented as a difference $g-h$ of two nondecreasing functions
$g, h$. A standard construction of such decomposition, i.e.,
$g(t)$ defined as a variation of $f$ on $[0,t]$, and $h=g-f$,
shows that we may additionally assume that $g$ maps $[0,1]$ into
$[0,1]$ and the image of $h$ is contained in $[-1,1]$. Therefore,
the image of the continuous map $\vf: H\times H\to
[-1,2]^{[0,1]}$, defined by $\vf(h_1,h_2) = h_1 - 2h_2 + 1$,
contains $BV$. Hence the space $BV$ is a continuous image of a
closed subset $K$ of the product $H\times H$. Clearly, there exist
measures $\mu\in P(BV)$ which are not strongly countably
determined, since the space $BV$ is not first countable. By
\cite[Proposition 2]{Po82}, every Radon probability measure on $K$
is strongly countably determined, and the above example shows that
this property is not preserved by continuous images.
\end{remark}

We shall now present our proof of a result due to Bourgain mentioned in the introductory section.

\begin{thm}\label{4:3}
Suppose that $X$ is a separable metrizable space and $K\sub \er^X$
is such a compact space that for every $g\in K$ and every closed
set $F\sub X$, the restricted function $g_{|F}$ has a point of
continuity. Then every measure $\mu\in P(K)$ is of countable type.
\end{thm}

\begin{mypr}
Let  us fix a measure $\mu\in P(K)$; for any $x\in X$ and $r\in\er$ we write
\[ L(x,r)=\{g\in K: g(x)\le r\}.\] 

\noindent {\bf (1)} We start by the following elementary
observation: suppose that $A,B\sub K$ are measurable set such that
$\mu(A\btu B)\ge\delta$ and $|\mu(A)-\mu(B)|\le \delta/2$; then
\[\mu(A\sm B)=\mu(A)-\mu(A\cap B)\ge \mu(B)-\mu(A\cap B)-\delta/2=\mu(B\sm A)-\delta/2,\mbox{ so}\]
\[2\cdot\mu(A\sm B)\ge \mu(A\btu B)-\delta/2\ge \delta -\delta/2=\delta/2,\]
and therefore $\mu(A\sm B)\ge\delta/4$.
\medskip

\noindent {\bf (2) Claim.} For a fixed $q\in\er$ there is a countable $X(q)\sub X$ such that
for every $x\in X$
\[\inf\{\mu(L(y,q)\btu L(x,q)): y\in X(q)\}=0.\]

Otherwise there is $\delta>0$ and an uncountable set $T\sub X$ such that
\[ (\dagger)\quad \mu( L(t,q)\btu L(s,q)))\ge\delta \mbox{ for } t,s\in T, t\neq s.\]
Shrinking $T$ if necessary, we can additionally assume that for $t,s\in T$
\[ (\ddagger)\quad |\mu(L(t,q))-\mu(L(s,q))|\le \delta/2.\]
For every $t\in T$ let us write $L_t^0=L(t,q)$ and choose $q_t>q$
so that $\mu(L_t^0)+\mu(L_t^1)\ge 1 - \delta/12$, where
\[L_t^1=\{g\in K: g(t)\ge q_t\}.\]
Note that, again choosing a suitable  uncountable subset of $T$, we can in fact
assume that $q_t=q'$ for every $t\in T$ and some fixed $q'>q$.

Now for every $t,s\in T$, if $t\neq s$ then
\[ \mu(L_t^0\cap L_s^1)\ge \mu(L_t^0\sm L_s^0)-\delta/12\ge \delta/4-\delta/12=\delta/6,\]
where we used (1) together with  $(\dagger)$ and  $(\ddagger)$.
This means that that we can apply Corollary \ref{3:4} with $\eps=\delta/6$; therefore
there is a dense-in-itself set $S\sub T$ such that
the family $\{(L_s^0,L_s^1)): s\in S\}$ is independent.

Let $P=\ol{S}$; we can divide $S$ into disjoint
subsets $S_0,S_1$ so that $\ol{S_0}=\ol{S_1}=P$. But then
\[C=\bigcap_{s\in S_0} L_s^0\cap\bigcap_{s\in S_1} L_s^1\neq\emptyset,\]
by independence and compactness. If we consider $g\in C$ then
$g(s)\le q$ for every $s\in S_0$, while $g(s)\ge q'>q$ for every $s\in S_1$,
so the function $g_{|P}$ has no point of continuity, a contradiction.
\medskip

\noindent {\bf (3)} We apply (2) to every $q\in Q$ and put $X_0=\bigcup_{q\in\qu} X(q)$.
It is routine to check that the countable algebra of sets generated by
$L(x,q)$, $x\in X_0$, $q\in\qu$ is $\btu$-dense in $Bor(K)$, and the proof is complete.
\end{mypr}

\begin{mysec}
Appendix: Proof of Theorem \ref{3:3}\label{appendix}
\end{mysec}

Let us note first that the assumptions of Theorem \ref{3:3} as well as the assertion of the result
can be expressed in terms of the measure algebra of $\mu$. If we denote by $\Sigma$ the $\sigma$-algebra
of subsets of $K$ generated by all $L_\xi^i$, $\xi<\omega_1$, $i=0,1$ then $\mu_{|\Sigma}$ is a measure
of type $\omega_1$ so, by the Maharam theorem, the corresponding measure algebra can be embedded
into the measure algebra of the usual product measure $\lambda$ on $2^{\omega_1}$.
Therefore it is sufficient to consider the case
when all $L_\xi^i$ are measurable subsets of $2^{\omega_1}$ and

\begin{liro}
\item $\lambda(L_\xi^0)+\lambda(L_\xi^1)>1-\eps/2$ for every $\xi<\omega_1$;
\item $ \lambda(L^0_\xi\cap L^1_\eta)> \eps$ whenever $\xi,\eta<\omega_1$, $\xi\neq \eta$.
\end{liro}

Let us recall that to prove Theorem \ref{3:3} we need to find an uncountable
$T\sub\omega_1$, and define a $ccc$-family $\AAA_0\sub \fin(T)$, such that
$\bigcup\AAA_0=T$ and $\{(L_\xi^0,L_\xi^1):\xi\in a\}$ is $\lambda$-independent whenever
$a\in\AAA_0$.

For  $I\sub \omega_1$, a set $B\sub 2^{\omega_1}$ is said to {\em be determined  by coordinates
in $I$} if $B=\pi_I^{-1}[\pi_I[B]]$, where $\pi_I:2^{\omega_1}\to 2^I$ denotes the projection; we
write $B\sim I$ to denote such a property.
Recall that we may think that $\lambda$ is defined on the Baire $\sigma$-algebra
$\algb=Baire(2^{\omega_1})$ of $2^{\omega_1}$,
which consists of sets of the form $B=B'\times 2^{\kappa\sm I}$, where $I\sub \omega_1$ is
countable and $B'\in Bor(2^I)$.
For any set $J\sub \omega_1$ we shall write $\algb[J]$ for the $\sigma$-algebra
of those $B\in\algb$ which are determined by coordinates in $J$.

As in \cite{Fr97}, for a set $B\in\algb$ and $J\sub\omega_1$ we denote by
$S_J(B)$ the set
\[S_J(B)=\bigcap_{I\sub J} (\chi_I\oplus B) =
2^{\omega_1}\sm \pi_{\omega_1\sm J}^{-1}[\pi_{\omega_1\sm
J}[2^{\omega_1}\sm B]],\] where $\oplus$ denotes the
coordinatewise addition mod 2. Note that  the set $S_J(B)$
is measurable and determined by coordinates in $\omega_1\sm J$.

Before we  start the main argument
we shall mention the following two auxiliary facts --- the first
one can be found in \cite{Fr97}.

\begin{lem}\label{5:0}
If the sets $J_n$ are
pairwise disjoint, $k\in\omega$, and, for every $n$, $|J_n|\le k$, then
\[\lim_n\lambda(S_{J_n}(B))=\lambda(B),\]
for every $B\in \algb$.
\end{lem}

\begin{lem}\label{5:1}
Let $B^0, B^1\in\algb$ be disjoint,  and let $J\sub\kappa$ be such that the set
\[\pi_J^{-1}\pi_J[B^0]\cap \pi_J^{-1}\pi_J[B^1],\]
has positive measure.
Then there is a finite set $c\sub\kappa\sm J$,
nonempty disjoint clopen sets $V^0,V^1\sim c$, and
a set $Z\in \algb[\kappa\sm c]$ with $\lambda(Z)>0$, such that
$Z\cap V^i\sub B^i$ for $i=0,1$.
\end{lem}

\begin{mypr}
If $\vf:c\to 2$ is a function defined on a finite set
$c\sub\omega_1$ then we write
\[C(\vf)=\{x\in 2^{\omega_1}: x_{|c}=\vf\},\]
for the cylinder set defined by $\vf$. Let us note that the
measure $\lambda$ satisfy the Lebesgue density theorem: \ if
$A\in\algb$ then for $\lambda$-almost all $x\in A$, $x$ is a
density point of $A$, i.e.\  we have
\[ \lim_{a}\lambda(A\cap C(x_{|a}))/\lambda(C(x_{|a}))=1,\]
where the limit operation is applied to a net directed by all finite sets (to see this use
the Lebesgue density theorem for the Cantor set and the fact that every $A\in\algb$ is determined by
countably many coordinates).

Our assumption on the sets $B^0,B^1$ implies that there are $x_i\in B^i$, $i=0,1$, such that
$\pi_J(x_0)=\pi_J(x_1)$, and in fact we can additionally assume that $x_i$ is a density point
of $B^i$, for $i=0,1$. From this we can conclude that there are clopen nonempty cylinders $W,V^0,V^1$
such that $W\sim J$, $V^0,V^1\sim \omega_1\sm J$, $V^0\cap V^1=\emptyset$, and for $i=0,1$
\[(\dagger)\quad \lambda(W\cap V^i\cap B^i)>(1/2)\lambda(W\cap V^i).\]
Moreover, we can take $V^i$ so that $V^i=C(\vf_i)$, where $\vf_0,\vf_1$ are defined on the same
finite set $c\sub\omega_1\sm J$.

Consider a function $f:2^{\omega_1}\to 2^{\omega_1}$ defined by
$f(x)=x\oplus (\vf_0\chi_c\oplus \vf_1\chi_c)$, where $\vf_i\chi_c$ denote the elements of
$2^{\omega_1}$ that extend $\vf_i$ by putting 0 outside $c$.
Such a function $f$ preserves the measure and $f[V^0]=V^1$,
so from $(\dagger)$ we get that the set
\[H=W\cap f[V^0\cap B^0]\cap V^1\cap B^1,\]
has positive measure; let
\[Z=\pi^{-1}_{\omega_1\sm c}\pi_{\omega_1\sm c}[H];\]
then $\lambda(Z)>0$ and $Z\sim \omega_1\sm c$. Now it suffices to check that
$Z\cap V^0\sub B^0$ and $Z\cap V^1\sub B^1$.

Let $x\in Z\cap V^0$; as $x\in Z$, there is $h\in H$ such that $h$ agrees with $x$ outside $c$.
In turn $h=f(y)$, where $y\in V^0\cap B^0$; $x$ agrees with $y$ outside $c$, while
$x_{|c}=\vf^0=y_{|c}$.  Finally, $x=y\in B^0$. We can check the
other inclusion in a similar way.
\end{mypr}

Using Lemma \ref{5:1} and our assumptions (i)-(ii) we construct inductively
an uncountable set $T\sub\omega_1$ and
\[ c_\xi\in \fin(\omega_1),\quad Z_\xi,V^0_\xi, V_\xi^1\sub 2^{\omega_1},\]
for $\xi\in T$ and $i=0,1$, so that the following are satisfied:

\begin{liar}
\item \label{1}
$\{c_\xi: \xi\in T\}$ is a pairwise disjoint family in $\fin(\omega_1)$;
\item \label{2}
$V_\xi^0, V_\xi^1$ are disjoint nonempty clopen sets and $V_\xi^i\sim c_\xi$ for $i=0,1$;
\item \label{3}
$\lambda(Z_\xi)>0$ and $Z_\xi\sim I_\xi$ for some countable $I_\xi\sub \omega_1\sm c_\xi$;
\item \label{4}
$V_\xi^i\cap Z_\xi\sub L_\xi^i$ for $i=0,1$.
\end{liar}

The inductive step can be done by the following observation: suppose that $S\sub \omega_1$
is a countable set and we have carried out the construction for $\eta\in S$. Then
we let
\[J=\bigcup_{\eta\in S} (c_\eta\cup I_\eta);\]
as $J$ is countable, $\lambda$ on $\algb[J]$ is of countable type, so there
must be $\xi\in\omega_1\sm S$ such that the sets $B^i=L^i_\xi$ satisfy the assumption
of Lemma \ref{5:1}.

Now we let $\AAA_0$ be the family of those finite sets $a\sub T$, for which
there is a set $Z\in\algb$ with $\lambda(Z)>0$, determined by coordinates
in a countable set $I\sub\omega_1$, so that
\[ I\cap\bigcup_{\xi\in a} c_\xi=\emptyset,\quad \mbox{and }
Z\sub \bigcap_{\xi\in a} Z_\xi.\]
Note that if $a\in\AAA_0$ and $\vf:a\to 2$ then,
\[\lambda\left(\bigcap_{\xi\in a} L_\xi^{\vf(\xi)}\right)\ge
\lambda\left(\bigcap_{\xi\in a} Z_\xi\cap V_\xi^{\vf(\xi)}\right)\ge
 \lambda\left(Z\cap \bigcap_{\xi\in a} V_\xi^{\vf(\xi)}\right)=
\lambda(Z)\cdot\prod_{\xi\in a}\lambda(V_\xi^{\vf(\xi)})>0,\]
where $Z$ is a set witnessing that $a\in\AAA_0$. This means that $\AAA_0$ consists
of sets $a$ making $\{(L_\xi^0,L_\xi^1): \xi\in a\}$ $\lambda$-independent;
now the proof of Theorem \ref{3:3} we be completed by the following fact which proof closely
follows Fremlin \cite{Fr97}, (see part (ii) of the proof of Theorem 6).

\begin{lem}\label{5:2}
$\AAA_0$ is a $ccc$-family.
\end{lem}

\begin{mypr}
Let $\{a_\beta: \beta<\omega_1\}\sub \AAA_0$ be a pairwise disjoint family; we can assume
that all the sets
\[ d_\beta=\bigcup_{\xi\in a_\beta} c_\xi,\]
are of constant size $k$.

Further we can assume that for every $\beta$, we have chosen a set
$Z(\beta)$, $Z(\beta)\sim I(\beta)$, witnessing that
$a_\beta\in\AAA_0$ so that for $\beta,\beta'<\omega_1$ we have
$\lambda(Z(\beta)\cap Z(\beta'))>\delta$, where $\delta>0$ is
fixed. This can be done using 2-linkedness of
measure algebras, see Lemma 6.16 in \cite{CN} (with  $n=2$).

Once we have done all those reductions,  there is $\beta>\omega$ such
that
\[d_\beta\cap \bigcup_{n<\omega} I(n)=\emptyset,\]
because $I(n)$ are countable and $d_\beta$ are pairwise disjoint.
By Lemma \ref{5:0} there is $n<\omega$ such that $\lambda(S_{d_n}(Z(\beta)))>\lambda(Z(\beta)))-\delta$, which
gives a nonnull set $W=Z(n)\cap S_{d_n}(Z(\beta))$. As $W$ is determined by coordinates
in $\omega_1\sm (d_n\cup d_\beta)$,  it follows
that $a_n\cup a_\beta\in\AAA_0$, and we are done. \end{mypr}

{\it Institute of Mathematics, University of Warsaw, ul.\  Banacha 2,  02--097 Warszawa,
Poland}

{\sc wmarcisz@mimuw.edu.pl}
\medskip

{\it Mathematical Institute, University of Wroc\l aw, pl.\ Grunwaldzki 2/4,
50-384 Wroc\-\l aw, Poland}

{\sc grzes@math.uni.wroc.pl} \quad
\verb#http://www.math.uni.wroc.pl/~grzes #

\end{document}